\documentclass[12pt]{article}
\usepackage{amsmath,amsfonts,amssymb,amsthm}
\usepackage{enumerate}
\usepackage{thm-restate}
\usepackage{fullpage}
\usepackage{cite}
\usepackage[margin=1in]{geometry}
\usepackage{tikz}
\usetikzlibrary{cd}

\theoremstyle{definition}
\newtheorem{theorem}{Theorem}[section]
\newtheorem{proposition}[theorem]{Proposition}
\newtheorem{corollary}[theorem]{Corollary}
\newtheorem{lemma}[theorem]{Lemma}
\theoremstyle{definition}
\newtheorem{definition}[theorem]{Definition}
\newtheorem*{remark}{Remark}
\newtheorem*{notation}{Notation}

\newtheorem*{example}{Example}

\newcommand{\N}{\mathbb N}
\renewcommand{\P}{\mathbb P}
\newcommand{\Q}{\mathbb Q}
\newcommand{\Z}{\mathbb Z}

\newcommand{\Spec}{\operatorname{Spec}}

\newcommand{\piet}{\pi_1^{\acute et}}

\newcommand{\inv}{^{-1}}

\newcommand{\OK}{\mathcal O_K}

\title{A Dynamical Analogue of the Criterion of N\'eron-Ogg-Shafarevich}
\author{Mark O.-S. Sing}

\begin{document}
\maketitle

\begin{abstract}
	We introduce an anabelian approach to the study of arboreal Galois representations and apply Tamagawa's anabelian version of the N\'eron-Ogg-Shafarevich criterion~\cite{tamagawa} to produce a dynamical analogue of this criterion: unramified representations correspond to rational maps satisfying a strong form of good reduction in terms of their critical locus. Subsequently, we pursue a dynamical anlaogue of the N\'eron-Ogg-Shafarevich criterion in terms of the more (dynamically) traditional arboreal representations, which relates unramified arboreal representations to a certain separability condition on the dynamical system. Finally, we relate the our criteria: the anabelian criterion corresponds to the dynamical criterion as one varies the base point around the critical locus. Along the way we develop effective criteria to determine which primes are infinitely ramified in arboreal representations over number fields, as well as the asymptotic growth of that ramification; we conclude with examples and applications, especially to dynamical systems over number fields.
\end{abstract}

\section{Introduction.}

The present paper has its origins in a combination of two different projects: the search for a dynamical analogue of the Neron-Ogg-Shafarevich (NOS) criterion~\cite{SerreTate} for the good reduction of abelian varieties, and extending the author's previous work on higher ramification in dynamical field extensions~\cite{DSen}. Both were accomplished by the author elementary methods under certain relatively mild restrictions (see preprint~\cite{Sing.DNOS.preprint}). Subsequently, the author realized that these results echoed work of Tamagawa~\cite{tamagawa} in anabelian geometry, and specifically Tamagawa's anabelian analogue of the criterion of N\'eron-Ogg-Shafarevich. In the present paper we elaborate on this connection: we begin by reformulating Tamagawa's anabelian criterion in the setting of dynamical systems, then connect this to the arboreal representations typically studied in arithmetic dynamics. This anabelian approach and the results we obtain suggest that, from the perspective of Galois representations, the traditional notion of good reduction in arithmetic dynamics is not quite strong enough, and that an alternative notion in terms of the pre-critical locus (or equivalently, how ``far'' the rational map is from being \'etale) may be more appropriate.

Ramification has always played an important role in the study of arboreal representations, beginning with the introduction of arboreal representations by Odoni~\cite{Odoni}. For example, progress towards of various versions of Odoni's Conjecture~\cite{BeneJuul,Looper,Specter} make essential use of ramification-theoretic methods. More recently, ramification in arboreal extensions has been studied in its own right~\cite{AHPW,DSen,JonesHamblen,bigSilverman,AitkenHajirMaire,BergerIterated,SpecterPHT}, with interesting applications to dynamical systems over both global and local fields.

In the analogy between dynamics and abelian varieties, arboreal representations are viewed as dynamical Tate modules. A source of the difficulty in studying arboreal representations is that they lack a counterpart for the Tate module's algebraic structure. With this in mind, Tamagawa's result is quite striking from the perspective of a dynamicist interested in arboreal representations: just like a dynamical system, a typical curve does not have any sort of algebraic structure, and so Tamagawa uses the pro-$\ell$ \'etale fundamental group as a replacement. This is exactly what we will endeavor to do in the the first part of this paper. Since we are concerned with dynamical systems on the curve $\P^1$, we will compare our dynamical systems to elliptic curves specifically. Our guiding analogies are as follows:

\begin{center}
\begin{tabular}{c|c}
	Elliptic Curves & Dynamical Systems on $\P^1$\\\hline
	      & $D$, the pre-critical pro-divisor\\
	 $U = E$  & $U = \P^1 - D$\\
	 $[\ell]:U \rightarrow U$&  $f:U\rightarrow U$\\
	 $[\ell]^*:\piet(U_{\bar K})\rightarrow \piet(U_{\bar K})$ & $\phi:f^*: \piet(U_{\bar K})\rightarrow \piet(U_{\bar K})$\\
	 $T_\ell = \piet(U_{\bar K})^{(\ell)}$ & $\piet(U_{\bar K})^{(d!')}$
\end{tabular}
\end{center}

In other words, we view $f$ and the multiplication-by-$\ell$ maps as \'etale endomorphisms of some variety, and use the geometric \'etale fundamental groups equipped with the induced endomorphisms as our ``Tate module''. For example, when $f$ is a Latt\'es map associated to a quotient $\pi:E\rightarrow \P^1$ from an elliptic curve $E$ and multiplication-by-$\ell$ map, the pre-critical divisor of $f$ is the image of $E[\ell^{\infty}]$ by $\pi$.

The reader should be careful to note that the arboreal representations with arbitrary base points \textit{do not} correspond to the Tate module in our framework: it is only the pre-critical representation, an intrinsic invariant of the dynamical system, which corresponds to the Tate module. This is not entirely surprising - for abelian varieties, the base point of the Tate module is always the identity. However, more general arboreal representations can still be incorporated into our anabelian framework. For example, it turns out that the pre-critical arboreal representation exhibits (almost) all of the possibilities for ramification in any arboreal representation -- a kind of ramification-theoretic upper bound.

To help motivate our results, let us recall a (\textit{very}) loose proof of the N\'eron-Ogg-Shafarevich criterion for elliptic curves and isolate its dynamical content.

\begin{center}
	$E$ has good reduction at $p$
	
	\hfill \ \ \ \ $\Longleftrightarrow$\hfill ($*$)
	
	the reduction of $E[\ell^n]$ modulo $p$ has $\ell^{2n}$ distinct elements
	
	\hfill \ \ \ \ $\Longleftrightarrow$\hfill ($**$)
		
	the Galois action on the $\ell$-adic Tate module is unramified
\end{center}

The dynamics we consider -- rational functions on $\P^1$ -- does not take place on a highly structured geometric object like the elliptic curve, and it is not necessarily clear what the dynamical analogue of an elliptic curve should be. As such, it is not clear that the equivalence ($*$) has a purely dynamical interpretation. However, the equivalence ($**$) amounts to a separability requirement for $[\ell^n]$: the preimage of $O$ should have as many elements mod $p$ as the degree of $\ell^n$. This separabilty requirement naturally extends to the dynamical setting: when does the preimage of $\alpha \in \bar K$ by $f^n$ have $\deg f^n$ distinct preimages in the residue field? When this happens, it follows immediately from Hensel's lemma that the associated arboreal representation is unramified, and our interest is in the converse. This is the primary concern of Section 3.

At the same time, this kind of separability condition plays a prominent role in Tamagawa's anabelian analogue of the N\'eron-Ogg-Shafarevich theorem. In Tamagawa's theory, good reduction of a curve $X$ punctured at a divisor $D$ requires an extension of the pair $(X,D)$ over $K$ to a pair $(\mathfrak X,\mathfrak D)$ over $\OK$ where $\mathfrak X$ is smooth and $\mathfrak D$ is relatively \'etale, which recovers the original pair $(X,D)$ on the generic fiber. One can immediately lift $\P^1_K$ to a smooth model $\P^1_{\OK}$, and not difficult to check that the divisor $D$ extends to a relatively \'etale divsior $\mathfrak D$ if and only if it has no components with multiplicity, and no points in its support reduce to the same point in the residue field after making a change of coordinate so that $0,1,\infty$ lie in its support. Tamagawa shows that an affine punctured curve has good reduction in this sense if and only if the Galois action on the geometric \'etale fundamental group is unramified.

These observations suggest that the \'etale fundamental group may provide a bridge between Tate modules and arboreal representations, and so it is natural to ask what can be done to make the endomorphism $f$ \'etale, so that it induces an endomorphism of the \'etale fundamental group. To do so, it suffices to delete the pre-critical locus. In combination, this suggests a close link between good reduction of the pre-critical locus (in the sense of Tamagawa) and ramification. Making this rigorous is the primary goal of this paper.

Our main results are best summarized by the following diagram. The Galois group of the ground field $K$ is denoted by $\Gamma$, the prime $p$ is the residue characteristic, and $\ell\neq p$ is an auxiliary prime.
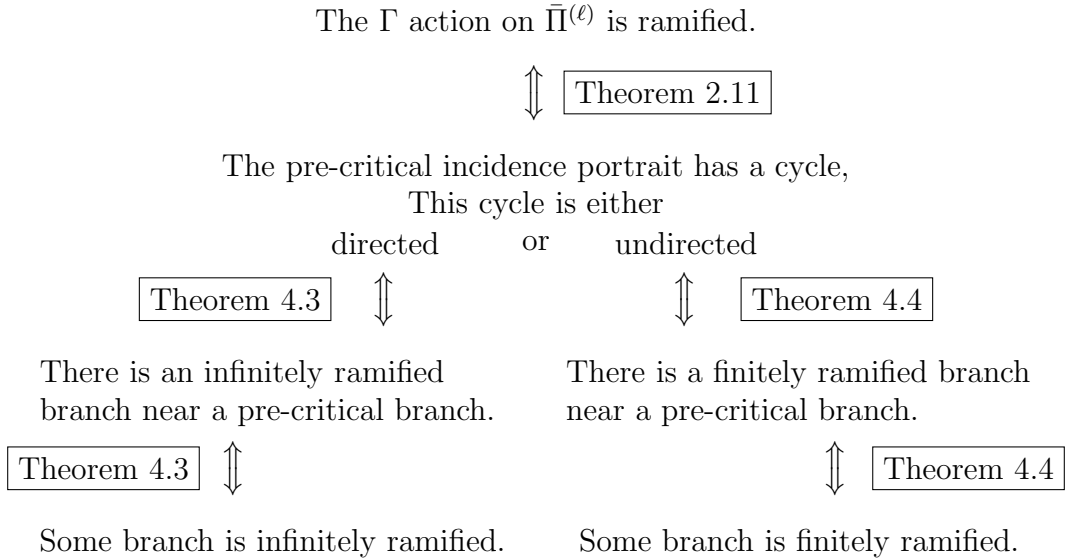
\begin{center}
\begin{tikzpicture}[scale=1]
	\draw node at (0,0) {The $\Gamma$ action on $\bar \Pi^{(\ell)}$ is ramified.};
	
	\draw node[rotate=90] at (0,-1) {$\Longleftrightarrow$};
	\draw node[rectangle,draw] at (1.75,-1) {Theorem~\ref{thm:nos1}};
	
	\draw node at (0,-2) {The pre-critical incidence portrait has a cycle,};
	\draw node at (0,-2.5) {This cycle is either};
	\draw node at (-2,-3) {directed};
	\draw node at (0,-3) {or};
	\draw node at (2,-3) {undirected};

	\draw node[rotate=90] at (-2,-3.75) {$\Longleftrightarrow$};
	\draw node[rectangle,draw] at (-4,-3.75) {Theorem~\ref{thm:main1}};
	
	\draw node[rotate=90] at (2,-3.75) {$\Longleftrightarrow$};
	\draw node[rectangle,draw] at (4,-3.75) {Theorem~\ref{thm:main2}};

	\draw node[text width=15em] at (-3.5,-5) {There is an infinitely ramified branch near a pre-critical branch.};
	\draw node[text width=15em,align=left] at (3.5,-5) {There is a finitely ramified branch near a pre-critical branch.};
	
	\draw node[rotate=90] at (-4,-6) {$\Longleftrightarrow$};
	\draw node[rectangle,draw] at (-5.75,-6) {Theorem~\ref{thm:main1}};
	\draw node at (-3.5,-7) {Some branch is infinitely ramified.};
	
	\draw node[rotate=90] at (4,-6) {$\Longleftrightarrow$};
	\draw node[rectangle,draw] at (5.75,-6) {Theorem~\ref{thm:main2}};
	\draw node at (3.5,-7) {Some branch is finitely ramified.};
\end{tikzpicture}
\end{center}

In Theorem 4.3, the infinitely ramified branch can be taken to be arbitrarily close to the pre-critical branch. If no critical points are periodic then one can take the infinitely ramified branch to be pre-critical. In contrast, for Theorem~\ref{thm:main2} the finitely ramified branch \textit{cannot} necessarily be chosen to be arbitrarily close to the pre-critical branch, and in general lies on a p-adic annulus around the pre-critical branch.

\newpage
\subsection{Notation.}

\begin{enumerate}[\indent--]
	\item $p$ is a prime,
	\item $K$ is a finite extension of $\Q_p$,
	\item $v$ is the valuation on $K$, normalized so $v(K)=\Z$,
	\item $\OK$ is the ring of integers of $K$ and $\pi_K$ a uniformizer of $\OK$,
	\item $\tilde K = \OK/\pi_K$ is the residue field of $K$, of characteristic $p$ and order $q=p^r$,
	\item $\bar K$ is a fixed algebraic closure of $K$,
	\item $X$ is the projective line $\P^1_K$ over $K$,
	\item $f(x) \in \OK(x)$ is a rational map of degree $d<p$ with good reduction,
	\item $D$ is the pro-divisor (see Section 2) associated to the pre-critical locus of $f(x)$,
	\item $U=\P^1_K - D$ as a $\Spec K$-scheme.,
	\item $\bar x$ is a geometric fixed point of $f(x)$,
	\item $\Gamma = \piet(\Spec U,\bar x)$ is the absolute Galois group of $K$.
\end{enumerate}

Except when $f$ is a powering map (a case we will not consider), the divisor $D$ always contains at least $3$ points. Making a change of coordinate if necessary, assume that $0$, $1$, and $\infty$ are in the support of $D$.

Much of this paper revolves around lifting objects defined over $\Spec K$ to objects defined over $\Spec \OK$. To simplify the presentation, we adopt the convention of using Latin letters for objects defined over $\Spec K$, and gothic letters for those defined over $\Spec \OK$.

We simply define $\mathfrak X = \P^1_{\OK}$ to be the gluing of $\OK[x]$ and $\OK[y]$ by $x\leftrightarrow 1/y$. Every point $a$ of $\P^1_K$ lifts to a point $
\mathfrak a$of $\P^1_{\OK}$ (take a minimal polynomial with integral coefficients not all in the maximal ideal and homogenize its minimal polynomial). Since $f(x)$ has good reduction, it extends to an endomorphism $\mathfrak f$ of $\mathfrak X = \P^1_{\OK}$. Likewise, we can extend $D$ to a relative pro-divisor $\mathfrak D$ over $\Spec \OK$. It will still be the case that $\mathfrak D$ is the pre-critical divisor of $\mathfrak f(x)$. More importantly, if the pair $(X,D)$ has good reduction in the sense of Tamagawa~\cite{tamagawa}, then the pair $(\mathfrak X,\mathfrak D)$ witnesses this fact, since we changed coordinate to place $0,1,$ and $\infty$ in $D$.

Given a point $\mathfrak a \in \mathfrak X$ associated to $a\in X$, we can restrict it to a point $\tilde a \in \P^1_{\tilde K}$ of the special fiber, called the \textit{reduction} of $a$. Observe that when $a\in\P^1(K)$ is not integral, it lifts to the point $1/a$ on the open affine $\Spec\OK[y,1/y]$ of $\mathfrak X$, and its reduction is the point at infinity in $\P^1_{\tilde K}$.

\newpage

\section{Anabelian Preliminaries.}

In this section, we establish the notation and facts necessary to extend Tamagawa's characterization of good reduction for finitely punctured curves to infinitely punctured curves. We will work in slightly greater generality at the start of this section, allowing $X$ to be an arbitrary scheme until otherwise noted.

Dynamical considerations typically give rise to infinite diagrams of schemes, and then to infinite projective limits in the category of schemes. Such limits need not exist in general, but they will when the morphisms in the diagram are affine. So we restrict our attention to dynamical systems which in which we can take such limits
\begin{definition}\label{def:nearly.etale}
	A dynamical system $f:X\rightarrow X$ is called \textbf{nearly \'etale} when it satisfies the following property: if $D_n$ is the ramification locus of $f^n$ and $U_n = X - D_n$, then the inclusion of $U_n$ into $X$ is an affine morphism (hence the inclusion of $U_{n+1}$ into $U_n$ is also affine).
\end{definition}

It is worth noting that for sufficiently nice schemes $X$ with a finite endomorphism $f:X\rightarrow X$, the associated dynamical system will typically be nearly \'etale. This is true of all curves, and in fact any curve over $K$ is itself affine after making at least one puncture at a $K$-rational point. Indeed, it can be shown that in many cases the purity of the ramification locus implies that inclusion of the complement $U_n$ of the ramification locus of $f^n$ is affine~\cite{StacksPure}.

In this sense, being nearly \'etale is a strong dynamical purity condition.

\begin{proposition}\label{prop:etale.restriction}
	If $f:X\rightarrow X$ is finite and nearly \'etale, then the restriction $f:U_{n+1}\rightarrow U_{n}$ is finite \'etale.
\end{proposition}
\begin{proof}
	Consider an open affine $\Spec A$ in $U_n$. Because $f$ and the inclusions of $U_n$ and $U_{n+1}$ in $X$ are affine, the inverse image of $\Spec A$ by $f$ is an affine open $\Spec B$ in $X$. Moreover, $f\inv (U_n) = U_{n+1}$ and so $\Spec B$ is an open affine in $U_{n+1}$. The extension of rings $B/A$ induced by $f$ is finite \'etale, so we are done.
\end{proof}

Dynamics on curves is essentially limited to genus zero and genus one; the latter is the domain of elliptic curves, so we are only concerned with the situation in genus zero. In this case, every finite endomorphism is nearly \'etale: 
\begin{proposition}\label{prop:P1.nearly.etale}
	Let $A$ be an integral domain and $f:\P_A^1 \rightarrow \P_A^1$ a dynamical system defined by a nonconstant rational map $f$ over $A$. This dynamical system is a nearly \'etale dynamical system. 
\end{proposition}
\begin{proof}
	Morphisms associated to nonconstant rational maps are finite, and upon dehomogenizing it is easy to see that the inclusions of each $D_n$ into $X$ are affine.
\end{proof}

Just as dynamical considerations require us to consider infinite diagrams of schemes, we will also want to consider divisors (arising, for instance, from orbits or backward orbits) which have infinitely many points. This requires a small generalization from divisors to pro-divisors, which essentially allow infinite formal sums with possibly infinite coefficients:
\begin{definition}\label{def:pro.divisor}
	Let $X$ be a scheme. An \textbf{effective pro-divisor} $D$ on $X$ is a sequence of effective Cartier divisors $D_n$ such that for all $m\leq n$, the scheme-theoretic intersection of $D_n$ and $D_m$ is $D_m$.
\end{definition}

A general philosophy in dynamics is that the critical orbit controls the behavior of the dynamical system; we are working ``backwards'' and so we work with the pre-critical orbit :
\begin{definition}\label{def:precrit.divisor}
	The \textbf{pre-critical pro-divisor of $f$}, denoted $C_f$, is the effective pro-divisor on $X$ determined by the system of ramification divisors $R_{f^n}$ of $f^n$. Note that $C_f = C_{f^n}$ for all $n$.
\end{definition}

The complement of $C_f$ ``should be'' the largest open subscheme on which all iterates of $f^n$ are \'etale. Of course, this complement is neither open nor a scheme in general. However, the complement of the the pre-critical divisor of a nearly \'etale endomorphism is a scheme, though not open.
\begin{theorem}\label{thm:lim.exists}
	Given a nearly \'etale dynamical system $f:X\rightarrow X$, the associated limit $U_\infty = \varprojlim U_n$ of the complements of the ramification divisors with their natural inclusions exists in the category of schemes, and $f:U_\infty\rightarrow U_\infty$ is finite \'etale.
\end{theorem}
\begin{proof}
	The inclusion of $U_n$ into $X$ is affine, hence also its inclusion into $U_{n-1}$. Thus the morphisms in the diagram are all affine, and in this case such a limit does exist in the category of schemes. The restrictions $f:U_n\rightarrow U_{n-1}$ are all finite \'etale, and so they have a limit $f:U_\infty\rightarrow U_\infty$ as well.
\end{proof}

\begin{proposition}\label{prop:cover.descent}
	If $Y\rightarrow U_\infty$ is a finite \'etale cover, then there is a sequence of finite \'etale covers $Y_n\rightarrow U_n$ with affine morphisms $Y_n\rightarrow Y_{n-1}$ such that $Y_n\rightarrow U_n$ is finite \'etale for all sufficiently large $n$, and $\varprojlim Y_n = Y$.
\end{proposition}
\begin{proof}
	We may assume that $Y$, $U_\infty$, and all $U_n$ are affine. Select finitely many generators $\{y_i\}$ for $\Gamma(Y,\mathcal O_Y)$ over $\Gamma(U_\infty,\mathcal O_{U_\infty})$ and adjoin them to $\Gamma(U_n, \mathcal O_{U_n})$; the spectrum of this ring is $Y_n$. The coefficients of the polynomial equations defining the $y_i$ are contained in large enough $U_n$, and their ramification divisors must be contained in the support of the pro-divisor $C_f$, hence in the ramification divisor of $f^n$ for all $n$ sufficiently large.
\end{proof}

\begin{theorem}~\label{thm:homotopy.exact}
	Suppose that $X$ is defined over a field $K$. Let $\bar x$ be a geometric fixed point of $f$ which is not pre-critical. Then 
	$$\pi_1(U_\infty,\bar x) = \varprojlim \pi_1(U_n,\bar x)$$
	and hence has an endomorphism $\phi:\pi_1(U_\infty,\bar x) \rightarrow \pi_1(U_\infty,\bar x)$ induced by $f$.
	
	Moreover, there is an \'etale homotopy exact sequence
\begin{center}
	\begin{tikzcd}
		0 \arrow[r]& (\pi_1(U_\infty)_{\bar K},\bar x) \arrow[r] & \pi_1(U_\infty,\bar x) \arrow[r] & \pi_1(\Spec K,\bar x) \arrow[r] & 0 
	\end{tikzcd}
\end{center}
	which is $\phi$-equivariant.
\end{theorem}
\begin{proof}
	The existence of these groups and the endomorphism $\phi$ follows from Theorem~\ref{thm:lim.exists}. The exact sequence arises as the limit of the \'etale homotopy exact sequences associated to the $U_n$. 
\end{proof}

When $f$ is a rational endomorphism of $\P^1_K$ and not a powering map, the pre-critical locus is infinite, and so the geometric fundamental group $\bar\Pi$ is the free profinite group on countably many generators. The generators correspond to inertia groups over the pre-critical divisor, and the action of $\phi$ is determined by its action on those points.

The exact sequence of Theorem~\ref{thm:homotopy.exact} induces an (outer) action of $\Gamma$ on $\bar \Pi$. One can recover the pre-critical arboreal representation associated to $f$ from this action: this action descends to an action on conjugacy class of inertia subgroups over each preimage $\alpha$, such that $\sigma\in \Gamma$ takes the conjugacy class over $\alpha$ to the conjugacy class over $\alpha^\sigma$.

Following Tamagawa, we make the following definition of good reduction for pro-divisors:
\begin{definition}
	Let $K$ be a local field with ring of integers $\mathcal O_K$. Let $f:X\rightarrow X$ be a nearly \'etale dynamical system over $K$ with pre-critical pro-divisor $C_f$. Then we say that the dynamical system has \textbf{good critical reduction} if and only if there is a smooth proper variety $\mathfrak X$ over $\OK$ and a relatively \'etale divisor $\mathfrak D$ such that on the generic fiber we recover $X$ and $C_f$.
\end{definition}

We will only work with $X=\P^1_K$, and in this case the extension to $\mathfrak X$ is unique up a change of coordinate over $\OK$ after moving three points to $\{0,1,\infty\}$. In this setting an effective divisor on $\mathfrak X$ is relatively  \'etale if and only if it (1) has no components with multiplicity, and (2) no components which meet on the special fiber. Condition (2) originates from the stipulation that a relatively \'etale divisor must be a relative normal crossings divisor.

Good critical reduction is a stronger notion than typical good reduction:
\begin{proposition}\label{def.critnormal}
	If $f:\P^1_K\rightarrow \P^1_K$ has good critical reduction, then it has good reduction in the usual sense.
\end{proposition}
\begin{proof}
	Suppose $f$ has bad reduction. Then there is a coordinate in which we can write $f(x) = \frac{q(x)}{p(x)}$ where $p,q$ have integral coefficients with at least one coefficient a unit, and $q$ and $p$ have a common root $\tilde c$ on the special fiber that is the reduction of a root $c$ of $q(x)$ in $K^{nr}$. Necessarily, $c$ is a critical point of $f$, and $f^2$ has two critical points $c_1, c_2$ which both reduces to preimages of $c$ by $f$ on the special fiber, and so $f^2$, hence also $f$, cannot have good critical reduction.
\end{proof}

Now we can reformulate Tamagawa's anabelian criterion for good reduction~\cite{tamagawa} into an anabelian criterion for the good critical reduction of dynamical systems on $\P^1$:
\begin{theorem}\label{thm:nos1}
	Let $K$ be a $p$-adic field with ring of integers $\mathcal O_K$, residue characteristic $p$, and $f:\P^1_K\rightarrow \P^1_K$ a nearly \'etale dynamical system. Then this dynamical system has good critical reduction if and only if the (outer) action of $\Gamma$ on $\pi_1((U_\infty)_{\bar K}, \bar x)^{(\ell)}$, the pro-$\ell$ completion of $\pi_1((U_\infty)_{\bar K}, \bar x))$ for some (all) $\ell\neq p$, is unramified.
\end{theorem}
\begin{proof}
	Extended to pro-divisors, Tamagawa's anabelian criterion for good reduction~\cite[5.3]{tamagawa} tells us that the Galois action on the pro-$\ell$ geometric fundamental group is unramified if and only if the dynamical system has good critical reduction, with smooth model $(\P^1_{\OK},\mathfrak D)$. Since the dynamical system has good critical reduction, $f$ has good reduction by Proposition~\ref{def.critnormal} and hence extends to an endomorphism $\mathfrak f:\mathfrak X\rightarrow\mathfrak X$. Since the map has good critical reduction, the pre-critical pro-divisor $C_{\mathfrak f}$ is relatively \'etale and its restriction to the generic fiber is $C_f$.
\end{proof}

In our situation, there is no notion of ``potentially good critical reduction'', as it is not possible to separate points meeting on the special fiber even after taking finite extensions. This is witnessed by the fact that if the action of $\Gamma$ on $\bar \Pi$ is ramified, then it is infinitely ramified. This is in contrast to arboreal representations, even pre-critical, which can be potentially unramified.

Readers familiar with Tamagawa's result are likely aware that it is straightforward to prove by elementary means in the genus $0$ case. In some sense this is the starting point for the next sections, but we track more dynamical data along the way. So, while not strictly necessary, we chose to reframe our results in this anabelian setting because (1) the common language makes it easier to compare the two kinds of representation, and better illuminates the geometric aspects of the results, and (2) this paper serves as a test case in anticipation of future work which makes more essential use of these tools.

It is also worth noting that Nakamura's anabelian weight filtration and rigidity results~\ref{nakamura} over number fields extend to this setting. While we are already mixing anabelian and dynamical ideas, it is particularly interesting to note that a \textit{dynamical} rigidity result is an important ingredient in Nakamura's rigidity result.

\section{Dynamical Preliminaries.}

In this section we characterize the ramification of ``branch extensions'' in terms of the power series expansion of $f($ near the base point.

\begin{definition}\label{def:branch}
	A \textbf{branch} for $f$ over $K$ is a sequence $(\alpha_n)_{n\in\N}$ of elements of $\bar K$, such that $f(\alpha_{n+1}) = \alpha_n$. The \textbf{base point}, $\alpha_0$, must be in $K$.
\end{definition}

Given any branch, its Galois orbit naturally has the structure of a tree. The natural coordinate-wise action of $\Gamma_K$ on this tree gives rise to a kind of representation:
\begin{definition}\label{def:branch.rep}
	Given a branch $B = (\alpha_n)_{n\in\N}$ for $f$, let $T$ be the tree formed by the Galois orbit of $B$. The associated \textbf{branch representation} is the homomorphism from $\Gamma_K$ to $\textrm{Aut}(T)$ induced by the action of $\Gamma$ on the branches.
\end{definition}

The full arboreal representation can naturally be written in terms of the branch representations associated to each Galois orbit of branches.

As usual, such representations are called unramified if they factor through the unramified quotient of $\Gamma_K$, and infinitely ramified if there is no finite extension $L$ of the base field $K$ such that the restriction of the representation to $\Gamma_L$ is unramified.

\begin{definition}\label{def:ram.index}
	Suppose a branch $B = (\alpha_n)_{n\in\N}$ is residually periodic of exact period $m$. Define the \textbf{ramification index of the branch}, denoted $e_B$, as the Weierstrass degree of $f^m(x)$ expanded as a power series around $x=\alpha_0$.
\end{definition}
\begin{definition}\label{def:height}
	The \textbf{height} of $f$ as the largest integer $h$ such that (residually) we can write $\tilde f(x) = \tilde Q(x^{p^h})$ for some rational function $Q(x)$. A rational function has positive height if and only the reduction of its derivative is identically zero.
\end{definition}
\begin{notation}
For simplicity, we fix a choice of $Q_f$ for $f$ such that $\tilde f(x) = \tilde Q_f(x^{p^h})$, where $h$ is the height of $f$.
\end{notation}

The significance of height zero is the following fact, immediate from the definitions:
\begin{proposition}
	If $f$ has height zero, then every critical point of its reduction $\tilde f$ is the reduction of a critical point of $f$.
\end{proposition}

The ramification index of a residually fixed branch is constant on the open disk of (spherical) radius $1$ around the branch. Our use of the word ``height'' is in reference to the height of an endomorphism of formal group. The reductions of rational maps with positive height exhibit some behavioral similarities with positive height endomorphisms of formal groups, in the spirit of the analogy between arithmetic dynamics and abelian vavarieties.

The following proposition is a generalization of~\cite[Proposition 2.3]{DSen}. We only apply it to the power series expansions of rational maps, but the more general case fits naturally into the study of dynamical systems on the open $p$-adic disk initiated by Lubin~\cite{LubinDynamical,LubinFlow}.

\begin{proposition}\label{prop:ramstable}
	Let $f(x)$ be a power series in $\OK[[x]]$ with finite Weierstrass degree $e$ and such that $f(0) = 0$. Let $(\alpha_n)_{n\in \N}$ be any nontrivial branch for $f(x)$ contained in the open unit disk. Then for all $n$ sufficiently large,
	\begin{enumerate}[\indent(a)]
		\item $v(\alpha_{n+1}) = \dfrac{v(\alpha_n)}{e}$,
		\item the sequence $(e^nv(\alpha_n))_{n\in \N}$ is eventually constant,
		\item $K(\alpha_{n+1})/K(\alpha_n)$ is totally ramified of degree $e$.
	\end{enumerate}
\end{proposition}
\begin{proof}
	The proof is identical to that of~\cite[Proposition 2.3]{DSen}, using Weierstrass degree in place of degree.
	
	In brief, if the branch is nontrivial then it has an entry with finite valuation. Examining the Newton polygon, we see that $v(\alpha_n)$ is strictly decreasing, and when it is sufficiently small (i.e. index sufficiently large) we are guaranteed $v(\alpha_{n+1}) = v(\alpha_n)/e$, at which point the sequence $(e^nv(\alpha_n))$ will be constant. This verifies (a) and (b).
	
	As for (c), it follows from (a) or (b): if the valuations grow in that fashion, then the ramification index of $K(\alpha_{n+1})/K(\alpha_n)$ must eventually always be $e$. 
\end{proof}

A related result, under the assumption $p\nmid e$, is proven by Ingram~\cite{Uniformization}, who exploits the fact that polynomials have a totally ramified fixed point at $\infty$. After making a change of coordinate to move $\infty$ to $0$, the resulting rational function can be written as a power series with Weierstrass degree equal to its degree, to which Proposition~\ref{prop:ramstable} applies. When $p\neq e$, Ingram shows that $f$ is analytically conjugate to a power map. His construction can be directly adapted to power series whose Weierstrass degree is not divisible by $p$.

Evidently, Proposition~\ref{prop:ramstable} will allow us to determine when a branch representation is infinitely ramified, and to describe the asymptotic growth of the ramification index if we can find coordinates where the \textit{rational} function $f$ admits such a power series expansion. It remains to develop some machinery to relate the ramification of a branch representation to the existence of such a coordinate.

\begin{definition}\label{def:stepwise.simple.red}
	Take a branch $(\alpha_n)_{n\in \N}$ for a rational map of height $0$. We say that the branch has \textbf{stepwise simple reduction} if, for all $n$, $\tilde \alpha_{n+1}$ is a simple root of $\tilde f(x) - \tilde \alpha_{n}$. Additionally, we say that a branch has \textbf{eventually stepwise simple reduction} if it has stepwise simple reduction after removing an initial segment.
\end{definition}

Evidently, stepwise simple reduction is separability condition, and hence is closely related to the reduction of the critical points of $f$.

\begin{lemma}\label{lem:insep.periodic}
	Assume $f$ has height zero, and let $(\alpha_n)_{n\in \N}$ be a branch which does not have stepwise simple reduction. Then there is a critical point $c$ of $f$ whose reduction $\tilde c$ is periodic, and the entire residual branch is residually periodic, with its entries given by the orbit of $\tilde c$.
\end{lemma}
\begin{proof}
	If the branch does not have stepwise simple reduction, then $\tilde \alpha_n$ is a multiple root of $\tilde f(x) - \tilde \alpha_{n-1}$ infinitely often. This can occur only when $\tilde \alpha_n$ a critical point of $\tilde f$.
	
	Hence for each $n$ there is a residual critical point $\tilde c_n$ of $f$ such that $\tilde c_n = \tilde \alpha_n$. Since $f$ has height zero, its derivative is nonzero and therefore $f$ has only finitely many residual critical points, and each is the reduction of a critical point of $f$. The sequence $(\tilde c_n)_{n\in \N}$ repeats one of those values infinitely often, and therefore that residual critical point is residually periodic, and the branch below any entry where it appears is periodic. Since the value reappears arbitrarily high in the residual branch, the whole residual branch is periodic.
\end{proof}

When the height $h$ is larger than zero, the derivative of $\tilde f(x) = \tilde Q_f(x^{p^h})$ vanishes, so every point is residually critical, and hence no branch can have stepwise simple reduction. On the other hand, any choice of $Q_f$ has height zero and hence $Q_f$ can have branches with good critical reduction.

Observe that, residually, $\tilde f$ is the composition of a rational map of height zero with a power of the absolute Frobenius. So it is then natural think of $f$ as being a residual twist of $Q_f$ by $\Phi^h$. This suggests a natural untwisting process for branches for rational maps of positive height:

\begin{lemma}
\label{lem:twist}
	Let $f(x)$ be a rational map of height $h$. Let $\Phi$ be the absolute Frobenius automorphism of the residue field. Assume that the coefficients of $\tilde f(x)$ are fixed by $\Phi^h$. Then there is a correspondence between residual branches for $\tilde f(x)$ and branches for $\tilde Q_f(x)$ over the residue field:
	 $$(\tilde \alpha_n)_{n\in\N}
	 \ \ \ \Longleftrightarrow\ \ \ 
	 (\Phi^{-hn}(\tilde \alpha_n))_{n\in\N}$$
\end{lemma} 

\begin{proof}
Since $\Phi^h$ fixes the coefficients of $\tilde f(x) = \tilde Q_f(x^{p^h})$, it commutes with $\tilde Q_f$. Then 
$$\tilde\alpha_{n-1} = \tilde f(\tilde\alpha_n) = \tilde Q_f(\tilde\alpha_n^{p^h}) = \tilde Q_f(\Phi^h(\tilde\alpha_n)) = \Phi^h \tilde Q_f(\tilde\alpha_n)$$

Inductively,
$$\tilde \alpha_0 = \tilde f^n(\tilde \alpha_n) = \Phi^{nh} (\tilde Q_f(\alpha_n))$$

Therefore, $(\Phi^{-nh}(\tilde \alpha_n))_{n\in\N}$ forms a branch (over the residue field) for $\tilde Q_(x)$. Reversing this process turns a branch for $\tilde Q_f$ into a branch for $\tilde P(x)$.	
\end{proof}

Adopting the convention that $\Phi^0 = \textrm{Id}$, it is natural to incorporate the ``untwisting'' into Definition~\ref{def:stepwise.simple.red} and combine Lemmas~\ref{lem:insep.periodic}~and~\ref{lem:twist} to treat all heights at once:
\begin{definition}\label{def:rss_full}
	Let $f$ be a rational function of height $h$. A branch $(\tilde \alpha_n)_{n\in\N}$ for $P(x)$ is said to have good critical reduction if the corresponding ``untwisted'' residue branch $(\Phi^{-hn}(\tilde \alpha_n))_{n\in\N}$ for the height zero rational map $\tilde Q_f$ has good critical reduction in the sense of Definition~\ref{def:stepwise.simple.red}. 
\end{definition}
\begin{proposition}\label{prop:rss_full}
	Suppose $f$ is a rational map with height $h$, so that $\tilde f(x) = \tilde Q_f(x^{p^h})$ and let $(\alpha_n)_{n\in \N}$ be a branch for $f$. If the branch does not have eventually stepwise simple reduction, then it is residually periodic.
\end{proposition}
\begin{proof}
	If $h=0$, this is exactly the statement of Lemma~\ref{lem:insep.periodic}, so assume $h$ is positive.

	Consider an iterate $f^k$, which will have height $hk$, so that residually $\tilde f^k(x) = \tilde R(x^{p^{hk}})$ where $R'(x)\neq 0$. This iterate has coefficients in the same ground field as $f$,  so it is possible to choose a $k$ such that $\Phi^{hk}$ fixes all the coefficients of $f$. 
	
	If the branch $(\Phi^{-hn}\tilde \alpha_n)_{n\in\N}$ for $\tilde Q_f$ is not residually stepwise simple for $\tilde Q_f$, then neither is the sub-branch $(\Phi^{-hnk}\tilde\alpha_{nk})_{n\in\N}$ for $\tilde S(x)$. Therefore the untwisted branch is residually periodic; in other words, contained in a finite extension of the residue field. The automorphism $\Phi^h$ of this residue field has finite order, and so the twisted branch is still periodic, though possibly of larger period.
	
	Since the subbranch is residually periodic, so too is the original branch; again the period may be larger.
\end{proof}

\section{Main Results.}

Our initial dynamical analogue of the N\'eron-Ogg-Shafarevich criterion, Theorem~\ref{thm:nos1} is essentially a direct adaptation of Tamagawa's anabelian result. There, consideration of the \'etale fundamental group suggests a natural form of good reduction not previously considered. In this section we will explore its connection with the more traditional arboreal and branch representations of arithmetic dynamics.

\begin{definition}\label{def:div.graph}
	Let $G$ be a reduced pro-divisor on $\P^1_K$, with support containing $\{0,1,\infty\}$. Then $G$ admits a natural extension to a pro-divisor $\mathfrak G$ on $\P^1_{\OK}$, and by reduction a pro-divisor $\tilde G$ on $\P^1_{\tilde K}$. We will enrich all three to with the structure of a directed multigraph. To $\mathfrak G$, we add two kinds of edges:
	\begin{enumerate}
		\item a directed edge $\mathfrak g\rightarrow\mathfrak  h$ if $\mathfrak f(\mathfrak g) = \mathfrak h$, and
		\item an undirected edge $\mathfrak g \leftrightarrow \mathfrak h$ if $\mathfrak g$ and $\mathfrak h$ meet on the special fiber.
	\end{enumerate}	
	
	The graphs $G$ and $\tilde G$ will only have vertices of the first type; edges $g\rightarrow h$ or $\tilde g\rightarrow \tilde h$ according to the action of $f$ or $\tilde f$. Observe that the vertices of $G$ can be identified with the subgraph of $\mathfrak G$ obtained by deleting all undirected edges, while $\tilde G$ can be obtained as a quotient of $\mathfrak G$ by identifying all vertices connected by undirected edges.
	
	We refer to $G$, $\mathfrak G$, and $\tilde G$ as the \textbf{portraits of $f$ on the divisor $G$}.
\end{definition}
\begin{definition}\label{def:dyn.graph}
	Let $f:\P^1\rightarrow \P^1$ be a nonconstant rational map, and denote by $C$, $\mathfrak C_{\mathfrak f}$, and $\tilde C_{\tilde f}$ the portraits of $f$ on its pre-critical pro-divisor. We will refer to these, particularly $\mathfrak C_{\mathfrak f}$, as the \textbf{dynamical incidence graph(s) of $f$} or the \textbf{dynamical portrait(s) of $f$}.
\end{definition}

This graph is clearly related to our anabelian construction: the graph $\mathfrak G$ has an undirected edge if and only if the representation on the pro-$\ell$ \'etale fundamental group of $\P^1-G$ is ramified. Moreover, if a branch fails to have stepwise simple reduction, then dynamical incidence graph for the associated arboreal representation has an undirected edge. We will show that ramification of the arboreal representation can be described in terms of this incidence graph.

\begin{theorem}\label{thm:main1}
	The following are equivalent:
	\begin{enumerate}[(a)]
		\item The dynamical incidence graph $\mathfrak C$ on the pre-critical divisor has a directed cycle (i.e. a cycle with at least one directed edge).
		\item For all $\epsilon > 0$, there an $\alpha \in K$ and critical point $c$ of $f$ such that $|\alpha-c| < \epsilon$ and a branch over $\alpha$ is infinitely ramified.
		\item Some branch representation associated to $f:X\rightarrow X$ is infinitely ramified.
	\end{enumerate}

	When (b) or (c) above occurs, the aforementioned branch $(\alpha_n)_{n\in\N}$ is residually periodic and there is a pre-critical branch $\{\gamma_n\}_{n\in\N}$ such that $|\alpha_n - \gamma_n|<1$ for all $n$. Let $m$ be the exact period of the reduction of the branch and $e$ the ramification index of $f$ on the branch (Definition~\ref{def:ram.index}. Then for all sufficiently large $n$, the extension $K(\alpha_{n+m})$ over $K(\alpha_n)$ is totally ramified of degree $e$. In other words, the ramification index of $K(\alpha_n)$ over $K$ grows like $Ce^{n/m}$ for some constant $C$.
	
	Moreover, if the critical point $c$ in part (b) is not periodic then one may take $\alpha = c$ as a base point.
\end{theorem}

\begin{proof}

(a)$\Rightarrow$(b) Suppose that $\mathfrak C$ has a directed cycle. Replacing $f$ by an iterate, we may assume that this cycle has a single directed edge -- in other words, it gives rise to a fixed point on the special fiber. Since we are working with the pre-critical locus, we may apply $f$ to push the cycle down until it contains a critical point of $f$. Without loss of generality we may assume that $f(0) = 0$ and that $\tilde 0$ is the fixed point (i.e. our cycle) on the special fiber. If $c$ is not periodic, let $\alpha=c$. Otherwise let $\alpha$ be any element of $K$ with sufficiently high valuation. This guarantees that no branch based at $\alpha$ is periodic.

	Now simply take preimages by $f$ while holding the special fiber fixed, to obtain a branch $(\alpha_n)_{n\in\N}$ which does not eventually have stepwise simple reduction, because $\tilde f(x) - \tilde\alpha_n = \tilde f(x)$ will have a multiple root at $\tilde 0$ for all $n$. Because the branch is not periodic it follows from Proposition~\ref{prop:ramstable}, applied to the power series expansion of $f$ around $0$, that the branch is infinitely ramified. In fact, the ramification index at each step eventually grows as a power of the ramification index of $f$ as a power series around $0$, and hence has infinite pro-$\ell$ component for some prime $\ell \leq \deg f < p$.

(b)$\Rightarrow$(c) Trivial.

(c)$\Rightarrow$(a) Suppose there is some infinitely ramified branch, not necessarily pre-critical. Then by Lemma~\ref{lem:insep.periodic}, the restriction of this branch to the special fiber is periodic and given by the orbit of the reduction of a critical point. This residually periodic critical orbit immediately gives rise to a directed cycle in the pre-critical incidence graph, as well as one in any arboreal representation containing the branch.
	\\
	
	Finally, in (b) and (c) we have just seen that the branch is residually periodic, and hence we can apply Proposition~\ref{prop:ramstable} to $f^m$, which verifeis the claim about the eventual ramification behavior of the branch.
\end{proof}
\begin{remark}
	One can view $\mathfrak f: \mathfrak X\rightarrow \mathfrak X$ as a family of dynamical systems over $\OK$. One can interpret Theorems~\ref{thm:nos1}~and~\ref{thm:main1} as the following correspondence between arithmetic and dynamics: $\mathfrak f$ is dynamically stable if and only if all arboreal representations in a neighborhood of the pre-critical locus are unramified, which occurs if and only if the representation of $\Gamma$ on $\piet(U,\bar x)^{(\ell)}$ is unramified.
	
	Over a global field, a rational map $f(x)\in K(x)$ can be viewed as a family of dynamical systems over the $S$-integers of $\OK$ for some finite set of primes $S$. It follows from work of Silverman~\cite{bigSilverman} that there are infinitely many primes not in $S$ modulo which a given critical point is periodic, and hence has infinitely ramified arboreal representation for a large family of base points. As such, $f$ cannot be dynamically stable as a family of rational maps over $\Spec \mathcal O_{K,S}$. This is somewhat reminiscent of the well-known fact that no abelian variety has good reduction over $\Spec \Z$, although our fact is true for all number fields and their $S$-integers. It appears that infinite ramification and bad (critical) reduction are characteristic of dynamical systems, in contrast to abelian varieties, which are ramified at only finitely many primes, and not necessarily infinitely ramified at those primes.
\end{remark}

Given an arboreal representation whose associated dynamical incidence graph has an undetected cycle, the representation on the \'etale fundamental group will be ramified, but the arboreal representation need not be. However, moving the base point allows this ramification to be detected, although it will only result in a finite amount of ramification.

\begin{theorem}\label{thm:main2}
	The following are equivalent:
	\begin{enumerate}[(a)]
		\item The pre-critical dynamical incidence graph has an undirected cycle.
		\item Some branch extension $(\alpha_n)_{n\in \N}$ for $f$ is ramified.
	\end{enumerate}

In fact, given any pre-critical branch $(\gamma_n)_{n\in\N}$ which has a member lying in an undirected cycle of the pre-critical incidence portrait, there is a ramified branch $(\alpha_n)_{n\in\N}$ such that $|\alpha_n-\gamma_n|< 1$ for all $n$.
\end{theorem}
\begin{proof} Without loss of generality, we may assume that there are no directed cycles on any of graphs in question, as then Theorem~\ref{thm:main1} immediately verifies the claim.

(a)$\Rightarrow$(b) 
	If there is an undirected cycle, then there is such a cycle with only two members, $\mathfrak q$ and $\mathfrak r$. We may assume $\mathfrak q$ and $\mathfrak r$ are preimages of the same critical point $\mathfrak c$. Then there are integers $m$ and $n$ such that $\mathfrak f^m(\mathfrak q),\mathfrak f^n(\mathfrak r) = \mathfrak c$. If $m\neq n$, this would give rise to a directed cycle, contrary to assumption, so $m=n$. Taking the minimal such $m$, we see that $\mathfrak f^{m-1}(\mathfrak q)$ and $\mathfrak f^{m-1}(\mathfrak r)$ are both preimages of $\mathfrak c$. Finally, after a possibly making a change of coordinates, assume $c=0$ and $q,r \neq \infty$.

	Now let $\alpha_0 \in K$ be any point such that $v_K(\alpha_0) + f(q)) = 1$. There is a branch for $f$ based at $\alpha_0$ which is ramified. Observe that the choice of $\alpha_0$ in combination with the fact that $f'(0) = 0$ guarantees that $f(x+q) - \alpha_0$ has a root with non-integer valuation (consider the Newton polygon after expanding $f(x + q)-\alpha_0$ as a power series around $0$). The corresponding root of $f(x) - \alpha_0$ therefore gives rise to a nontrivial ramified extension, and any branch containing it will suffice.

(b)$\Rightarrow$(a) If a branch $(\alpha_n)_{n\in \N}$ is ramified, then the reduction of some $f(x) - \alpha_n$ has a double root. This double root must be a residual critical point of $f(x)$, and hence gives rise to an undirected edge (hence undirected cycle) in the pre-critical dynamical incidence graph.
\\

The final remark follows from the fact that the branches in question are all residually pre-critical, and hence entrywise close to a pre-critical branch. However, the ramified branch cannot always be made arbitrarily to a pre-critical branch, as happens when the pre-critical branch is not ramified.

\end{proof}

Combined, Theorems~\ref{thm:main1}~and~\ref{thm:main2} tell us that, in a sense, ``all'' of the possibilities for ramification associated to branches for $f$ can be understood by looking at the pre-critical locus and branches based on annulus around it.

\begin{theorem}
	If an arboreal representation is infinitely ramified, one of its branch representations is infinitely ramified.
\end{theorem}
\begin{proof}
	Evidently if some branch is infinitely ramified, so too is the arboreal representation.
	
	Let $\mathfrak T$ be the dynamical incidence portrait associated to the preimage tree. If no branch representation is infinitely ramified, then $\mathfrak T$ can have no directed cycles by Theorem~\ref{thm:main1}. As such, a given branch can ramify at most $2d-2$ times, corresponding to the residual critical points lying on the branch. Therefore, the ramification index is bounded uniformly by $(d!)^{2d-2}$ across branches. This ramification is tame because we have assumed $p > d$, so by Abhyankar's lemma the ramification in all the branches can be trivialized by replacing $K$ with any totally tamely ramified extension of degree $(d!)^{2d-2}$. Thus the full arboreal extension, obtained as the compositum of all the branch extensions, is unramified after a finite base change, and therefore only finitely ramified over $K$.
\end{proof}
\begin{remark}
	The previous theorem holds when the ramification is wild. One must take care to show that after an unramified base change which is \textit{independent of the branch}, the some initial segment of the branch is totally ramified of bounded degree, and afterward unramified. A $p$-adic field has only finitely many extensions of bounded degree and so there is a single finite base change which trivializes all the possibilities for ramification.
\end{remark}

\section{Applications}

\subsection{Post-Critically Finite Maps}
Our results are especially powerful when applied to post-critically finite maps over number fields.

Combining the main result of this paper with that of \cite{DSen}, we can exactly and effectively determine which primes are infinitely ramified in branch representations of of post-critically finite polynomials of prime-power degree with good reduction.
\begin{corollary}
\label{cor:pcf}
	Suppose $f$ is a post-critically finite rational map with good reduction defined over a number field $K$. Let $\Delta_{f}$ be the (necessarily finite) set $\bigcup_{f'(c) = 0}\{c - f^n(c)\ |\ n\in \N_+\}$. Then a prime $\mathfrak p$ of $K$ for which $f$ has height zero is infinitely in ramified in some arboreal representation if and only if $\mathfrak p$ divides a member of $\Delta_f$.
	
	In fact, an arboreal representation with base point $\alpha \in K$ is infinitely ramified at a prime $\mathfrak p$ for which $f$ has height zero if and only if there is a critical point $c$ of $f$ which is periodic modulo $\mathfrak p$ and $\alpha$ is in the orbit of $c$ modulo $\mathfrak p$.
\end{corollary}

There are only finitely many primes for which $f$ has height zero, as such a prime has to divide its degree.

\begin{example}
	Let $c$ be a root of $t^3 + 2t^2 + t + 1$, and set $f(x) = x^2 + c$. This is the post-critically finite map associated to the Douady rabbit. The critical orbit is $$\{0,c,c^2+c,\infty\},$$ where the first three points form a periodic cycle and the last is fixed.
	
	This polynomial has everywhere good reduction, and only has positive height for primes lying over $2$. Since all the critical points are periodic, for any prime $\mathfrak p$ not dividing $2$ there is an arboreal representation for $f$ which is infinitely ramified at $p$. The base point for this representation can be chosen to be $\mathfrak p$-adically near any periodic point of $f$.
	
	For any particular base point $\alpha\in K$ which is not periodic for $f$, the arboreal representation can only ramify at one of the finitely many prime $\mathfrak p$ which divide one of $\alpha, \alpha - c, \alpha - c^2 - c$, or $\alpha - \infty$ where a prime $\mathfrak p$ divides the latter if and only if $\frac 1 \alpha$ is divisible by $\mathfrak p$.
	
	As far as ramification at divisors of $2$, this is a post-critically finite map of prime degree, so the main result of~\cite{DSen} implies all such primes are infinitely ramified, and in fact deeply wildly ramified in a precise way.
	
	To be even more concrete, if we take $\alpha = 5$ then it is straightforward to check that the arboreal representation is infinitely ramified at exactly the following primes: $2$, and the ideals $(5-a)$, $(2\alpha^2 + 3\alpha + 1)$, $(\alpha^2 + 2\alpha + 3)$, which lie over $181$, $7$ and $19$, respectively.
\end{example}

\subsection{Post-Critically Infinite Maps}

The pre-critical arboreal representations of post-critically infinite maps over global fields exhibit surprising behavior in contrast to abelian varieties: these arboreal representations are infinitely ramified at infinitely many primes, while the representation on the Tate module of an abelian variety can only ramify at finitely may primes.

\begin{proposition}\label{prop:pci}
	Let $f$ be a rational map defined over a number field $K$, and $c$ a critical point of $f$ whose orbit is infinite. Then the arboreal representation associated to $f$ and $c$ is ramified at infinitely many primes.
\end{proposition}
\begin{proof}
	It follows from results of Silverman~\cite{Silverman} that there are infinitely many primes modulo which $c$ is periodic. Restricting to primes of good reduction and where $f$ has height zero, we obtain infinitely many primes $\mathfrak q$ such that $\tilde c$ is periodic. The relationship $\tilde f^n (\tilde c) = \tilde c$ gives rise to a directed cycle on the pre-critical dynamical incidence graph, and since $c$ itself is not periodic, Theorem~\ref{thm:main1} guarantees that a branch over $c$ is infinitely ramified at $\mathfrak q$.
\end{proof}

However, it is the author's suspicion that, at least for ``most'' rational maps, only representations containing a pre-critical representation (in other words, arboreal representations attached to post-critical points) can be infinitely ramified at infinitely many primes. The reason is that for such primes, not only must there be a critical point which is periodic modulo infinitely many primes, but the residual orbits of that critical point must contain the base point in their orbits as well. In other words, there must be infinitely many primes which divide an entry from both 
$$\{P^n(c) - c\ |\ n\in \N\}\ \ \ \ \textrm{and}\ \ \ \ \{P^m(c)-\alpha_0\ |\ m\in\N\}$$
Both sets are fairly sparse, and explicit computer calculations seem to indicate that they are unlikely have many prime divisors in common.

\subsection{Higher Ramification}

The methods of this paper can, in some cases, allow us to extend our previous work on higher ramification~\cite{DSen} to rational functions. This is particularly interesting from the perspective of developing a dynamical analogue of the theory of crystalline representations, which naturally arises when trying to relate good reduction at $p$ to the representation on the $p$-adic Tate module, rather than an $\ell$-adic Tate module for $\ell\neq p$.

\begin{theorem}
	Suppose $f$ is a rational function with good reduction and a periodic point $a$ with period $m$ such that $f^m(x+a)-a$ has Weierstrass degree a prime power. Take a branch $(\alpha_n)_{n\in\N}$ in the open unit disk around $a$.
	
	Then the associated arboreal extension is infinitely wildly ramified, and there is a polynomial $W$ determined by $f$, $m$, and $a$, such that all results of~\cite{DSen} on the structure of the higher ramification groups apply to $f^m(x+a)-a$ when $W$ is used to calculate the limiting ramification data.
\end{theorem}
\begin{proof}
	It follows immediately from Theorem~\ref{thm:main1} that any nontrivial branch representation in this setting is infinitely wildly ramified.
	
	As to the final claim, replace $f$ by $f^m$ and change coordinate so that $a=0$. Expand $f$ as a power series around $0$ -- this series hash integral coefficients, and converges and equals $a$ on the open unit disk. Apply the $p$-adic Weierstrass preparation theorem to write this series as $W(x)u(x)$, where $W$ is a polynomial with good reduction of the form $x^{p^h}$ modulo $\pi_K$ and $u(0)$ is a unit.
	
	Then observe that in all arguments of~\cite{DSen}, the Newton (co)polygon which determines the limiting ramification data is $f(x+\alpha_n)$, and because $u(0)$ is a unit, this is the same as the Newton (co)polygon of $W(x+\alpha_n)$.	
\end{proof}

\subsection{Abelian Dynamical Extensions}

Abelian dynamical extensions, especially over global fields, have attracted significant interest recently~\cite{PetscheAndrews,FOZ}. In fact, it is known that if an arboreal representation for a rational map $f(x)$ over a number field is abelian, then $f(x)$ must be post-critically finite~\cite{FOZ}. We can recover this by purely local methods in the special case of pre-critical arboreal representations: if any critical point has infinite orbit, then Proposition~\ref{prop:pci} furnishes infinitely many primes at which the arboreal representation is infinitely tamely ramified. But no such extension can be abelian -- in fact, the associated decomposition subgroup already fails to be abelian, and even remains nonabelian after any finite base change. So not only does this arboreal representation fail to be abelian, it is quite farm from abelian: it has infinitely many infinite nonabelian decomposition groups.

Of the known examples of abelian arboreal representations, many are pre-critical, such as the Latt\'es maps appearing in explicit class field theory of totally imaginary quadratic fields, or the trivial case of a powering map based at zero.

Moreover, the constraints we obtain on ramification are particularly precise for PCF rational maps. In combination with global class field theory, this goes a long way towards restricting the possibilities for abelian arboreal representations.

\section*{Acknowledgments}
I would like to thank my advisor, Joseph H. Silverman, for his mentorship and guidance, and also Laura DeMarco for asking the questions which led me to the remark after Theorem~\ref{thm:main1}.

\bibliography{refs}{}
\bibliographystyle{plain}

\end{document}